\documentclass[]{interact}

\usepackage{epstopdf}
\usepackage{subfigure}

\usepackage{natbib}
\bibpunct[, ]{(}{)}{;}{a}{}{,}

\theoremstyle{plain}
\newtheorem{theorem}{Theorem}[section]
\newtheorem{lemma}[theorem]{Lemma}

\newtheorem{proposition}[theorem]{Proposition}

\theoremstyle{definition}

\theoremstyle{remark}

\theoremstyle{observation}

\usepackage{amsmath,amssymb}
\usepackage{url}
\usepackage[mathscr]{euscript}
\usepackage{stmaryrd}
\usepackage{mathabx}
\usepackage{dsfont}
\usepackage{graphicx}
\usepackage{multirow}
\usepackage{xcolor,colortbl}

\newcommand{\scrC}{\mathscr{C}}

\newcommand{\scrO}{\mathscr{O}}
\newcommand{\scrP}{\mathscr{P}}
\newcommand{\scrS}{\mathscr{S}}
\newcommand{\bbN}{\mathbb{N}}
\newcommand{\bbR}{\mathbb{R}}
\newcommand{\bbone}{\mathds{1}}

\newcommand{\Proj}{\operatorname{Proj}}
\newcommand\defeq{\mathrel{\overset{\makebox[0pt]{\mbox{\normalfont\tiny\sffamily def}}}{=}}}

\begin{document}


\title{Beating the SDP bound for the floor layout problem: A simple combinatorial idea}

\author{
\name{Joey Huchette\textsuperscript{a}\thanks{CONTACT J. Huchette. Email: huchette@mit.edu}, Santanu S. Dey\textsuperscript{b}, and Juan Pablo Vielma\textsuperscript{c}}
\affil{\textsuperscript{a}Operations Research Center, Massachusetts Institute of Technology, Cambridge, MA\\ \textsuperscript{b}School of Industrial and Systems Engineering, Georgia Tech, Atlanta, GA\\ \textsuperscript{c}Sloan School of Management, Massachusetts Institute of Technology, Cambridge, MA}
}

\maketitle

\begin{abstract}
For many mixed-integer programming (MIP) problems, high-quality dual bounds can be obtained either through advanced formulation techniques coupled with a state-of-the-art MIP solver, or through semidefinite programming (SDP) relaxation hierarchies. In this paper, we introduce an alternative bounding approach that exploits the ``combinatorial implosion'' effect by solving portions of the original problem and aggregating this information to obtain a global dual bound. We apply this technique to the one-dimensional and two-dimensional floor layout problems and compare it with the bounds generated by both state-of-the-art MIP solvers and by SDP relaxations. Specifically, we prove that the bounds obtained through the proposed technique are at least as good as those obtained through SDP relaxations, and present computational results that these bounds can be significantly stronger and easier to compute than these alternative strategies, particularly for very difficult problem instances.
\end{abstract}




%

\section{Introduction}
A fundamental concept in optimization is the dual bound, which provides a proof on the quality of a given primal solution. Matching dual bounds (for minimization, a lower bound on the optimal cost) and primal bounds (e.g. the objective cost of a feasible solution) immediately provide a certificate of optimality, and for difficult optimization problems, tight bounds can provide confidence that an available solution is of sufficiently high quality.

Advances in computational methods \cite{DBLP:journals/anor/BixbyR07} such as cutting-plane technology~\cite{marchand:ma:we:wo:2002,RichardDey} and formulation techniques \cite{Mixed-Integer-Linear-Programming-Formulation-Techniques} have made the generation of good bounds with a state-of-the-art mixed integer programming  (MIP) solver attainable for a great number of problems. However, there are still many classes of problems for which the generation of bounds remains a challenge. A common approach to generate bounds for such problems is to use semidefinite programming (SDP) techniques to construct hierarchies of relaxations that theoretically converge to the best possible bound (e.g. \cite{Laurent:2002} and \cite[Section 10]{Conforti:2014}). Such SDP relaxations produce high quality bounds, but in practice can lead to semidefinite optimization instances which are too large to build in memory on a computer, much less solve to optimality.

In this work, we present an alternative combinatorial technique for constructing dual bounds that can sometimes significantly outperform, with regards to time and quality, the bounds obtained by both state-of-the-art MIP solvers and SDP relaxations. The technique is applicable to some specially structured MIP problems and exploits the \emph{combinatorial implosion} effect present in MIP and combinatorial optimization problems: that is,  slightly reducing the size of a problem can significantly reduce its solution time. To exploit this, the technique solves a series of fixed-sized subproblems obtained by only considering portions of the original problem and then aggregating the information to produce the global bound. This gives a polynomial-time bounding scheme that can be trivially parallelized and compares favorably against existing bounding techniques for certain problems, both theoretically and practically. In particular, we prove that when the technique is applied to a floor layout problem (FLP), the bounds generated are equal or better than those obtained by SDP relaxations. We also present computational results that show that the technique can produce better bounds in less time than both SDP relaxations and state-of-the-art MIP solvers.


\section{Combinatorial dual bounds}

The dual bounding scheme attempts to take advantage of \emph{combinatorial implosion}; that is, given a hard optimization problem, a smaller instance is likely much easier to solve than a slightly larger one\footnote{Combinatorial implosion is a more optimistic corollary to the well-known \emph{combinatorial explosion} effect, which is used throughout the folklore to describe the difficulty of many discrete optimization problems.}. To exploit this phenomena, we will take a difficult optimization problem and decompose it into many smaller optimization problems, which we may solve in a decentralized manner. Then, re-aggregating these subproblems, we may produce a dual bound on our original problem. We will perform this decomposition along certain subsets of the variables, which gives the scheme a combinatorial flavor.

For example, consider an optimization problem of the form $\gamma^* \defeq \min_{x \in Q} c^Tx$ for some set $Q \subseteq \bbR^n$ and $c \in \bbR^n$.
For some $S \subset \llbracket n \rrbracket \defeq \{1,\ldots,n\}$ and $u \in \bbR^n$, let $u_S \in \bbR^{|S|}$ represent the projection of $u$ on the components corresponding to $S$. If the problem $\min_{x \in Q} c^Tx$ is particularly difficult, we could instead consider solving a restricted version of the problem over variables $S \subset \llbracket n \rrbracket$, given by $\gamma(S,P_S) \defeq \min_{x_S \in P_S} c_S^Tx_S$ (notationally, we take $\gamma(S) \defeq \gamma(S,\Proj_S(Q))$). Provided that $P_S \supseteq \Proj_S(Q) \defeq \{ x_S : x \in Q \}$, this gives us a lower bound on a portion of the objective function in the following way.

\begin{lemma} \label{lemma:single-set}
    If $S \subset \llbracket n \rrbracket$ and $P_S \supseteq \Proj_S(Q)$, then
    \[
        \gamma(S,P_S) \leq \min_{x \in Q} c_S^Tx_S \leq c_S^Tx_S^*,
    \]
    where $x^* \in \arg\min_{x \in Q} c^Tx$.
\end{lemma}
\proof{}
    The first inequality follows as $\Proj_S(Q) \subseteq P_S$; the second, as $x^* \in Q$.
\endproof

We can reinterpret Lemma \ref{lemma:single-set} as producing an inequality $\gamma(S,P_S) \leq c_S^Tx^*_S$ that is valid for any feasible point $x \in Q$. With this interpretation, we can imagine producing inequalities for many such sets $S \subset \llbracket n \rrbracket$ and aggregating them into an optimization problem which will provide a dual bound for our original problem.

\begin{theorem} \label{thm:combinatorial-bound}
    Consider an optimization problem of the form $\gamma^* = \min_{x \in Q} c^Tx$ for some set $Q \subseteq \bbR^n$, and assume that $\gamma^* > -\infty$. Take some family $\scrS \subseteq 2^{\llbracket n \rrbracket}$ and corresponding relaxation sets $P_S$ with $P_S \supseteq \Proj_S(Q)$ for each $S \in \scrS$. Then for any $R \supseteq Q$, $\gamma^* \geq \omega(\scrS,\{P_S\}_{S\in\scrS},R)$, where
    \begin{subequations} \label{eqn:combinatorial-bound}
    \begin{align}
        \omega(\scrS,\{P_S\}_{S\in\scrS},R) \defeq \min_{x \in R} \quad& c^Tx \\
                   \operatorname{s.t.} \quad& c_S^Tx_S \geq \gamma(S,P_S) \quad \forall S \in \scrS. \label{eqn:combinatorial-bound-constraints}
    \end{align}
    \end{subequations}
\end{theorem}
\proof{}
    If the original problem is infeasible ($Q = \emptyset$), then this statement is true by the convention that $\min_{x \in \emptyset} c^Tx = \infty$. If the problem is bounded, consider an optimal $x^* \in \arg\min_{x \in Q} c^Tx$; that is, $c^Tx^* = \gamma^*$. We see that, from construction, $x^*_S \in P_S$ for each $S \in \scrS$, and therefore $c^T_S x^*_S \geq \gamma(S,P_S)$ for all $S \in \scrS$. Therefore, $x^*$ is feasible for the relaxed problem, and the bound follows.
\endproof

In a slight abuse of notation, we will take $\omega(\scrS) \defeq \omega(\scrS,\{\Proj_S(Q)\}_{S \in \scrS},R)$ when $R$ is understood to be a fixed relaxation for $Q$ (that is, fixed for different choices of $\scrS$).

\subsection{Applicability of Theorem~\ref{thm:combinatorial-bound}}
The potential strength of the bound from \eqref{eqn:combinatorial-bound} depends on the variables sets $\scrS$, but also the sets $R$ and $P_S$ we choose to optimize over. At one extreme, if we choose $R=Q$, \eqref{eqn:combinatorial-bound} reduces to our original optimization problem $\min_{x \in Q} c^Tx$ with a new set of redundant inequalities \eqref{eqn:combinatorial-bound-constraints}, and so \eqref{eqn:combinatorial-bound} will produce a tight bound but will likely be at least as difficult to solve as the original problem. On the other extreme, if we take $R = \bbR^n$, then the strength of the bound is derived exclusively from the partial bounds \eqref{eqn:combinatorial-bound-constraints}. A natural choice if $Q$ is a mixed-integer set is to take $R$ as the \emph{relaxation} of $Q$ with integrality conditions dropped, which will generally be relatively easy to optimize over, while offering some additional strength to the bound \eqref{eqn:combinatorial-bound}.

Constructing the sets $P_S$ for $S \in \scrS$ also poses an interesting challenge. Clearly if we take $P_S = \bbR^S$, we have $\gamma(S,P_S) = -\infty$, and so the strength of \eqref{eqn:combinatorial-bound} will be solely determined by optimizing over the set $R$. Fortunately, for many optimization problems such as the floor layout problem, it is straightforward to construct sets $P_S$ by considering a family of optimization problems related to $Q$, but over smaller variable sets $S$. For example, for the FLP, we may construct $P_S$ by considering a smaller instance of the FLP with a smaller number of components, but otherwise the same problem data. In this way, we can avoid solving the original, difficult optimization problem, and compute a dual bound by instead solving many smaller instances of the same problem on a series of restricted index sets.

\subsection{Computational advantages}
The combinatorial bounding scheme from Theorem \ref{thm:combinatorial-bound} involves computing a series of partial dual bounds and aggregating them. The bounding values $\gamma(S,P_S)$ can be computed in whatever convenient way, but if the sets $P_S$ are MIP-representable, a natural approach is just to use an off-the-shelf MIP solver. These computations can be carried out in parallel, and the set $\scrS$ can be iteratively enlarged to produce tighter bounds as needed. Indeed, our experiments provide an instance where, if $R$ is simple enough to optimize over, the vast majority ($>99\%$) of the computational time is spent in computing the bounds $\gamma(S,P_S)$ which can be done completely in parallel before the relatively cheap master problem \eqref{eqn:combinatorial-bound}.


\subsection{Computational complexity}
Similar to lifted hierarchies such as the Sherali-Adams or Lasserre hierarchies, this combinatorial bounding scheme yields a polynomial-time dual bounding scheme for appropriate choices of $\scrS$. More precisely, consider the nested families $\scrS_k = \{S \subseteq \llbracket n \rrbracket : |S| \leq k\}$ and the case where $Q$ is a mixed-integer linear set. The corresponding combinatorial bound for level $k$ corresponds to solving $\scrO(n^k)$ subproblems, each of which is a mixed-integer problem with at most $k$ binary variables. Therefore, for fixed $k$, the computational complexity of the combinatorial bounding scheme grows polynomially in the dimension $n$. Furthermore, when $k=n$, we recover a tight bound on the optimal cost.

\begin{proposition}
    The combinatorial bounding scheme produces a hierarchy of dual bounds in the sense that, for any nested sets $\scrS_1 \subset \scrS_2 \subset \cdots \scrS_{t-1} \subset \scrS_t = 2^{\llbracket n \rrbracket}$,
    \[
        \omega(\scrS_1) \leq \omega(\scrS_2) \leq \cdots \leq \omega(\scrS_{t-1}) \leq \omega(\scrS_t) = \gamma^*.
    \]
\end{proposition}


\section{Floor layout problem}\label{FLP_sec}
The floor layout problem (FLP), sometimes also known as the facility layout problem, is a difficult design problem traditionally arising in the arrangement of factory floors and, more recently, the very-large-scale integration (VLSI) of computer chips. FLP consists of laying out $N$ rectangular components on a fixed rectangular floor $[0,L^x] \times [0,L^y]$ in such a way as to minimize the weighted sum of communication costs between the components. The components do not have fixed dimensions, but must have fixed area, differentiating this problem from a more traditional 2D packing problem.

A number of mixed-integer formulations for the FLP have been proposed in the literature \cite{Bazaraa:1975,Meller:1999,Meller:2007,Sherali:2003,Castillo:2005,Huchette:2015}, though much of the work on the problem has been devoted to tailored heuristics \cite{Camp:1991,Anjos:2006,Bernardi:2010,Bernardi:2013,Jankovits:2011,Lin:2011,Luo:2008a,Liu:2007}. These heuristic approaches are often able to produce high quality solutions for large instances of FLP that are beyond the scope of exact MIP approaches. However, technology to obtain strong dual bounds is not nearly as mature. To the best of our knowledge, \cite{Takouda:2005} presents the most recent work on producing bounds for the FLP.


We will work with the unary formulation \eqref{eqn:2D-FLP} of FLP from \cite{Huchette:2015}, which is closely related to the FLP2 formulation from \cite{Meller:1999}; we refer the reader to \cite{Huchette:2015} for a more thorough discussion of alternative formulations. Notationally, we take $C$ as the set of components (for the moment, $C = \llbracket N \rrbracket$) and $P(C)$ as the set of all pairs of components in $C$.
\begin{subequations} \label{eqn:2D-FLP}
    \begin{alignat}{2}
        \min_{c,d,\ell,z}\quad& \sum_{(i,j) \in P(C)} p_{i,j} (d^x_{i,j}+d^y_{i,j})  & \label{eqn:2D-FLP-1} \\
        \operatorname{s.t.} \quad& d^s_{i,j} \geq c^s_i - c^s_j \quad& \forall s \in \{x,y\}, (i,j) \in P(C) \label{eqn:2D-FLP-2} \\
                         & d^s_{i,j} \geq c^s_j - c^s_i \quad& \forall s \in \{x,y\}, (i,j) \in P(C) \label{eqn:2D-FLP-3} \\
                         & \ell^s_i \leq c^s_i \leq L-\ell^s_i \quad& \forall s \in \{x,y\}, i \in C \label{eqn:2D-FLP-4} \\
                         & 4\ell^x_i \ell^y_i \geq \alpha_i \quad& \forall s \in \{x,y\}, i \in C \label{eqn:2D-FLP-5} \\
                         & c^s_i + \ell^s_i \leq c^s_j - \ell^s_j + L(1-z^s_{i,j}) \quad& \forall s \in \{x,y\}, i,j \in C, i \neq j \label{eqn:2D-FLP-6} \\
                         & z^x_{i,j} + z^x_{j,i} + z^y_{i,j} + z^y_{j,i} = 1 \quad& \forall (i,j) \in P(C) \label{eqn:2D-FLP-7} \\
                         & d^s_{i,j} \geq 0 \quad& \forall s \in \{x,y\}, (i,j) \in P(C) \label{eqn:2D-FLP-8} \\
                         & lb^s_i \leq \ell^s_i \leq ub^s_i \quad& \forall s \in \{x,y\}, i \in C \label{eqn:2D-FLP-9} \\
                         & 0 \leq z^s_{i,j} \leq 1 \quad& \forall s \in \{x,y\}, i,j \in C, i \neq j. \label{eqn:2D-FLP-10} \\
                         & z^s_{i,j} \in \{0,1\} \quad& \forall s \in \{x,y\}, i,j \in C, i \neq j. \label{eqn:2D-FLP-11}
    \end{alignat}
\end{subequations}

In words, formulation \eqref{eqn:2D-FLP} aims to optimally determine the centers $(c^x_i,c^y_j)$ and half-widths $(l^x_i,l^y_i)$ of components $i \in C$ such to minimize the weighted pairwise ``Manhattan'' norm distance $\sum_{(i,j) \in P(C)} p_{i,j}(|c^x_i-c^y_i| + |c^y_i-c^y_j|)$ (we assume that $p_{i,j} \geq 0$). Using the standard transformation, we linearize the absolute values in the objective 
by introducing auxiliary ``objective'' variables $d_{i,j}^s$, 
with constraints (\ref{eqn:2D-FLP-2}-\ref{eqn:2D-FLP-3}) that lower bound the corresponding absolute value term.
Constraint \eqref{eqn:2D-FLP-4} ensures each component lies completely on the floor, while \eqref{eqn:2D-FLP-5} enforces the area of each component is at least some constant $\alpha_{i}$. Together, (\ref{eqn:2D-FLP-6},\ref{eqn:2D-FLP-7},\ref{eqn:2D-FLP-11}) enforce that, for every pair of components, they are separated in (at least) one of the directions $x$ and $y$; this ensures that the resulting feasible solutions correspond to physical layouts, where none of the components overlap. 

For the remainder of the analysis, we assume the width of the floor is sufficiently large, in the sense that every possible layout is feasible; that is, $L^s \geq \sum_{i=1}^N ub^s_i$. Again, we remind the reader that $p$ is nonnegative, and so there is no incentive to place components farther apart.
%
%

\subsection{One-dimensional floor layout problem} The one-dimensional floor layout problem (1D-FLP) is the restriction of the FLP to a single direction: say, $x$. That is, it is the FLP where we impose that $z^y_{i,j}=0$ for all $(i,j) \in P(C)$. When restricted to a single direction, there always will be an optimal solution with $\ell^x_i = lb^x_i$ for all $i \in C$, so we remove the widths as decision variables from the 1D-FLP.
Equivalently, 1D-FLP is the problem of finding an optimal permutation of $N$ components in a line, such as to minimize the weighted sum of pairwise distances between the components. 

Besides being an interesting problem in its own right, the 1D-FLP
is a simpler problem that encompasses many of the challenges of solving the FLP.

A natural MIP formulation \eqref{eqn:1D-FLP} for the 1D-FLP, due to \cite{Love:1976},
is obtained by restricting the (two-dimensional) floor layout problem to a single dimension and fixing the widths to their lower bounds.
\begin{subequations} \label{eqn:1D-FLP}
    \begin{alignat}{2}
        \min_{c,d,z}\quad& \sum_{i,j} p_{i,j} d_{i,j}  & \label{eqn:1D-FLP-1} \\
        \operatorname{s.t.} \quad& d_{i,j} \geq c_i - c_j \quad& \forall (i,j) \in P(C) \label{eqn:1D-FLP-2} \\
                         & d_{i,j} \geq c_j - c_i \quad& \forall (i,j) \in P(C) \label{eqn:1D-FLP-3} \\
                         & c_i + \ell_i \leq c_j - \ell_j + L(1-z_{i,j}) \quad& \forall i,j \in C, i \neq j \label{eqn:1D-FLP-4} \\
                         & z_{i,j} + z_{j,i} = 1 \quad& \forall (i,j) \in P(C) \label{eqn:1D-FLP-5} \\
                         & 0 \leq c_i \leq L \quad& \forall i \in C \label{eqn:1D-FLP-6} \\
                         & d_{i,j} \geq 0 \quad& \forall (i,j) \in P(C) \label{eqn:1D-FLP-7} \\
                         & 0 \leq z_{i,j} \leq 1 \quad& \forall i,j \in C, i \neq j \label{eqn:1D-FLP-8} \\
                         & z_{i,j} \in \{0,1\} \quad& \forall i,j \in C, i \neq j. \label{eqn:1D-FLP-9}
    \end{alignat}
\end{subequations}

We note that a fairly extensive body of literature has studied the 1D-FLP, in terms of IP formulations \cite{Amaral:2006,Amaral:2008,Amaral:2012}, dual bounds \cite{Amaral:2009}, and SDP approaches \cite{Anjos:2005,Anjos:2008,Anjos:2009}. However, these formulations do not readily generalize to the two-dimensional FLP, and so we will maintain our attention in this work on formulations \eqref{eqn:2D-FLP} and \eqref{eqn:1D-FLP}.

\subsection{Combinatorial bounding scheme}\label{comb_sec}
We 
now specialize our dual bounding scheme for the 1D-FLP and FLP. For a given instance, consider the restricted problem on the components $C \subseteq \llbracket N \rrbracket$ and take the corresponding variables $S(C) \defeq \{(c_i, d_{j,k}, z_{j,k}, z_{k,j}) : i \in C, (j,k) \in \scrP(C)\}$ for the 1D-FLP and $S(C) \defeq \{(c_i, d_{j,k}, \ell_i, z_{j,k}, z_{k,j}) : i \in C, (j,k) \in \scrP(C)\}$ for the FLP; that is, all variables corresponding to components in the set $C$.

\begin{proposition}\label{prop:kbound}
    Consider some family $\scrC \subseteq 2^{\llbracket N \rrbracket}$. Then $\gamma^* \geq \omega(\scrC)$ for both the 1D-FLP and FLP, where
    \begin{subequations} \label{eqn:1D-FLP-bounding-scheme}
    \begin{align}
        \omega(\scrC) = \min_{d \geq 0} \quad& \sum_{(i,j) \in P(\llbracket N \rrbracket)} p_{i,j} d_{i,j} \\
            \operatorname{s.t.} \quad& \sum_{(i,j) \in P(C)} p_{i,j}d_{i,j} \geq \gamma(S(C)) \quad \forall C \in \scrC.
    \end{align}
    \end{subequations}
\end{proposition}
\proof{}
    First we consider the 1D-FLP. The result follows from Theorem \ref{thm:combinatorial-bound} by observing that $R = \bbR_+^C \times \bbR_+^{\scrP(C)} \times \bbR_+^{2\scrP(C)}$ outer approximates \eqref{eqn:1D-FLP}. Since the $c$ and $z$ variables do not appear in the objective, we may project them out and work solely on the objective variables $d$. For each $C \in \scrC$, we may take formulation \eqref{eqn:1D-FLP} over components $C$ as the outer approximation $P_S \supseteq \Proj_S(Q)$.

    We can apply the same argument to the FLP, but also may aggregate the variables $d_{i,j} \defeq d^x_{i,j} + d^y_{i,j}$.
\endproof

In particular, we will be interested in the nested families $\scrC_k \defeq \{C \subseteq \llbracket N \rrbracket : |C| \leq k \}$. These provide a natural hierarchy for studying this bounding scheme. We will use the notation $\omega_k \defeq \omega(\scrC_k)$. Furthermore, we note for the case of $\scrC_2$, the values $\gamma(S(\{i,j\})) = \ell_i+\ell_j$ for the 1D-FLP and $\gamma((S\{i,j\})) = \min\{\ell^x_i+\ell^x_j,\ell^y_i+\ell^y_j\}$ for the FLP, respectively, are available in closed form. Therefore, $\omega_{2}$ is also readily available as $\omega_2 = \sum_{(i,j) \in P(C)} p_{i,j}(\ell_i+\ell_j)$ or $\omega_2 = \sum_{(i,j) \in P(C)} p_{i,j}\min\{\ell^x_i+\ell^x_j,\ell^y_i+\ell^y_j\}$, respectively. We will be using this particular value for comparison later.

Finally, we note that, since many of the FLP instances instances have relatively sparse objectives (i.e. many $(i,j) \in \scrP(C)$ have $p_{i,j} = 0$), we may remove some sets from $\scrC_k$ without affecting the quality of the bound. In particular, if there is some $r \in C \in \scrC_k$ such that $p_{r,s} = 0$ for all $s \in C \backslash \{r\}$, then $\gamma(S(C)) = \gamma(S(C \backslash \{r\}))$ and also $\sum_{(i,j) \in \scrP} p_{i,j}d_{i,j} = \sum_{(i,j) \in \scrP(C\backslash\{r\})} p_{i,j}d_{i,j}$ for any feasible $d$. Since $C\setminus \{r\} \in \scrC_k$, we may omit $C$ from $\scrC_k$. In practice, this simple observation allows us to trim a significant number of sets, and is useful computationally.

\section{Comparison to other dual bounds}\label{comparizon_sec}
With our combinatorial bounding scheme for the 1D-FLP, we now wish to compare the approach to other standard approaches for producing dual bounds for mixed-integer optimization problems.



\subsection{LP relaxation quality} \label{ss:lp-relaxation}
It is not hard to see that the formulation \eqref{eqn:1D-FLP} produces a very weak dual bound; that is, the point
\begin{alignat*}{2}
    \hat{c}_i &\leftarrow L/2 \quad &\forall i \in C \\
    \hat{d}_{i,j} &\leftarrow 0 \quad &\forall (i,j) \in \scrP(C) \\
    \hat{z}_{i,j} &\leftarrow 1/2 \quad &\forall i,j \in C, i \neq j
\end{alignat*}
is feasible for the relaxation $\min\{\eqref{eqn:1D-FLP-1} : (\ref{eqn:1D-FLP-2}-\ref{eqn:1D-FLP-8})\}$ and has objective value zero. In other words, the LP relaxation provides no strengthening of the dual bound over the trivial observation that $\min_{d \geq 0}\sum_{(i,j) \in P(C)} p_{i,j}d_{i,j} = 0$ when $p_{i,j} \geq 0$.

We can improve the formulation by tightening the bounds on the objective variables:
\begin{equation} \label{eqn:obj-cut}
    d_{i,j} \geq \ell_i + \ell_j \quad \forall i,j \in C, i \neq j.
\end{equation}
With these inequalities, the LP relaxation value now matches with the combinatorial bound $\omega_2$.

Similarly, for the FLP we have that the optimal cost of the continuous relaxation $\min\{\eqref{eqn:2D-FLP-1} : (\ref{eqn:2D-FLP-2} - \ref{eqn:2D-FLP-10})\}$ of \eqref{eqn:2D-FLP} is zero, by considering the feasible solution
\begin{alignat*}{2}
    \hat{c}^s_i &\leftarrow L^s/2 \quad &\forall s \in \{x,y\}, i \in C \\
    \hat{d}^s_{i,j} &\leftarrow 0 \quad &\forall s \in \{x,y\}, (i,j) \in \scrP(C) \\
    \hat{\ell}^s_i &\leftarrow \max\{lb^s_i,\sqrt{\alpha_i}\} \quad &\forall s \in \{x,y\}, i \in C \\
    \hat{z}^s_{i,j} &\leftarrow 1/4 \quad &\forall s \in \{x,y\}, i,j \in C, i \neq j.
\end{alignat*}

Similarly, imposing the valid inequalities
\begin{equation} \label{eqn:obj-cut-2D-FLP}
    d^x_{i,j} + d^y_{i,j} \geq \min\{lb^x_i+lb^x_j,lb^y_i+lb^y_j\} \quad \forall (i,j) \in P(C)
\end{equation}
gives a relaxation bound equal to the first level of the combinatorial bounding scheme, $\omega_2$.

\subsection{Sparse valid inequalities for the 1D-FLP}
Based on the observation that \eqref{eqn:obj-cut} and \eqref{eqn:obj-cut-2D-FLP} help improve the bound of the LP relaxation, we might hope that we can identify some clever inequalities that produce stronger dual bounds. However, we show next that any set of valid inequalities (of support-wise disjoint families of valid inequalities; see details below) that significantly tighten the dual bounds for 1D-FLP in comparison to $\omega_k$ must be correspondingly dense. Sparse cutting-planes have good computational properties~\cite{Walter14} and are preferred by solvers, and recently there has been considerable work towards understanding the strength of sparse cutting-planes in general~\cite{deyMolinaroWang:2015,deyMolinaroIroume:2015,deyMolinaroWang:2016}.


Let $C \subseteq \llbracket N \rrbracket$. We call a valid (linear) inequality for \eqref{eqn:1D-FLP} as having \emph{support over $C$} if the only nonzero coefficients correspond to the decision variables $S(C)$.


\begin{proposition}\label{prop:validinequalities}
	Take $k,t \in \bbN$ such that $kt \leq N$. Let $C^1, C^2, \dots, C^t \subseteq \llbracket N \rrbracket$ such that $|C^l| = k$ for all $l \in \llbracket t \rrbracket$ and the $C^l$'s are pairwise disjoint. Consider the linear programming relaxation of formulation \eqref{eqn:1D-FLP} augmented with any number of valid inequalities where each valid inequality added has a support over some $C^l$ for $l \in \llbracket t \rrbracket$. Then the optimal cost of the resulting linear program is no greater than $\omega_{k}$.
\end{proposition}
\proof{}

    Assume w.l.o.g. that $C^1 = \{1,\ldots,k\}$, $C^2 = \{k +1, \ldots, {2k}\}$, $\dots$, $C^t = \{k(t-1) + 1, \dots, kt\}$. Moreover w.l.o.g. we may assume that an optimal permutation of the 1D-FLP over components $C^l$ is given by $(k(l-1) + 1,\ldots, kl)$.

Now consider two feasible solutions of 1D-FLP $(c^1, d^1, z^1)$ and $(c^2, d^2, z^2)$ given by
    \begin{alignat*}{2}
        c^1_i &= \ell_i + 2\sum_{j=1}^{i-1} \ell_j &\quad \forall i \in \llbracket N \rrbracket \\
        c^2_i &= \ell_i + 2\sum_{j=i+1}^N   \ell_j &\quad \forall i \in \llbracket N \rrbracket \\
        d^1_{i,j} &= |c^1_j - c^1_i| &\quad \forall (i,j) \in P(\llbracket N \rrbracket) \\
        d^2_{i,j} &= |c^2_j - c^2_i| &\quad \forall (i,j) \in P(\llbracket N \rrbracket) \\
        z^1_{i,j} &= \bbone[i < j] &\quad \forall i,j \in \llbracket N \rrbracket, i \neq j \\
        z^2_{i,j} &= \bbone[i < j] &\quad \forall i,j \in \llbracket N \rrbracket, i \neq j.
    \end{alignat*}
    Since these are feasible layouts, they must be feasible w.r.t. \eqref{eqn:1D-FLP} and the valid inequalities. Call the midpoint of these two points as $(\hat{c},\hat{d},\hat{z}) \defeq \frac{1}{2}(c^1,d^1,z^1) + \frac{1}{2}(c^2,d^2,z^2)$, i.e., $(\hat{c},\hat{d},\hat{z})$  is the point
    \begin{alignat*}{2}
        \hat{c}_i &= L/2 &\quad \forall i \in \llbracket N \rrbracket \\
        \hat{d}_{i,j} &= |c^1_j - c^1_i| &\quad \forall (i,j) \in P(\llbracket N \rrbracket) \\
        \hat{z}_{i,j} &= 1/2 &\quad \forall i,j \in \llbracket N \rrbracket, i \neq j.
    \end{alignat*}
    Since $(\hat{c},\hat{d},\hat{z})$ is a convex combination of two valid points, this point satisfies the valid inequalities that are added.

 Now examine a new point $(\tilde{c},\tilde{d},\tilde{z})$ where all the $d$-variables which do not belong to the support of $S(C^j)$ for some $j \in \llbracket t \rrbracket$ are set to $0$ as follows:
    \begin{alignat*}{2}
        \tilde{c}_i &= L/2 &\quad \forall i \in \llbracket N \rrbracket \\
        \tilde{d}_{i,j} &= \begin{cases} |c^1_j - c^1_i| & i,j \in C^l \textup{ for some } l\in \llbracket t \rrbracket \\ 0 & \text{o.w.} \end{cases} &\quad \forall (i,j) \in P(\llbracket N \rrbracket) \\
        \tilde{z}_{i,j} &= 1/2 &\quad \forall i,j \in \llbracket N \rrbracket, i \neq j.
    \end{alignat*}
    Note two properties enjoyed by  $(\tilde{c},\tilde{d},\tilde{z})$:
    \begin{enumerate}
        \item It satisfies all the constraints (except the integrality requirement on $z$-variables) in the formulation \eqref{eqn:1D-FLP}.
        \item It is indistinguishable from solution $(\hat{c},\hat{d},\hat{z})$ on the variables $\bigcup_{j \in \llbracket t \rrbracket}S(C^j)$ and therefore satisfies the valid inequalities that are added.
    \end{enumerate}
    Hence $(\tilde{c},\tilde{d},\tilde{z})$ is a valid solution to the linear program consisting of \eqref{eqn:1D-FLP} together with the valid inequalities. Finally, observe that by our choice of $(c^1, d^1, z^1)$ and $(c^2, d^2, z^2)$ with respect to the optimal permutation of 1D-FLP over the components $C$, we have that the objective function value of $(\tilde{c},\tilde{d},\tilde{z})$ is equal to $\sum_{l = 1}^t\gamma(S(C^l))$. Hence the optimal cost of \eqref{eqn:1D-FLP} augmented with the valid inequalities is at most $\sum_{l = 1}^t\gamma(S(C^l))$. It remains to show that $\gamma_k \geq \sum_{l = 1}^t\gamma(S(C^l))$ to complete the proof.

    Refer to Proposition \ref{prop:kbound} and let $d^{*}$ be an optimal solution of the LP corresponding to computation of $\omega_k$. Then note that
    \begin{eqnarray*}
    \gamma_k= \sum_{(i,j) \in P(\llbracket N \rrbracket)} p_{i,j}d^{*}_{i,j} &\geq& \sum_{l = 1}^t  \sum_{(i,j) \in P(C^l)} p_{i,j}d^{*}_{i,j} \geq     \sum_{l = 1}^t\gamma(S(C^l)),
    \end{eqnarray*}
    where the first inequality follows from non-negativity of $p_{i,j}$ and $d^{*}_{i,j}$ and the assumption that $C^l$'s are pairwise disjoint, and the last inequality follows from the constraints satisfied by $d^{*}$.
\endproof

Proposition \ref{prop:validinequalities} shows that, unless we use a large number of valid inequalities with non-disjoint supports, we are unlikely to beat the bound given by the combinatorial technique.

\subsection{Lifted relaxations for the 1D-FLP}
Next we consider the standard approach
of constructing a lifted relaxation with auxiliary variables whose projection onto the original space is (hopefully) a tight approximation of the feasible region. There are a host of frameworks for constructing these lifted representations in a generic way, with a natural distinction between those that result in linear optimization problems 
\cite{Balas:1993,Sherali:1990}) and semidefinite optimization problems \cite{Lovasz:1991,Lasserre:2001}.

The literature applying these techniques to combinatorial optimization problems---yielding both positive and negative results---is incredibly vast; we refer interested readers to the surveys of \cite{Chlamtac:2012} and \cite{Laurent:2002}, as well as Chapter 10 of \cite{Conforti:2014}.

We will adapt these techniques to the mixed-integer setting for the 1D-FLP. A natural construction is to consider $A\lambda\leq b$ as the system of linear inequalities defining the relaxation of \eqref{eqn:1D-FLP}, (with $\lambda = (c,d,z)$) and construct the nonlinear system
\begin{align*}
    z_{i,j}(A\lambda-b) &\geq 0 \quad \forall i,j \in \llbracket N \rrbracket, i \neq j \\
    (1-z_{i,j})(A\lambda-b) &\geq 0 \quad \forall i,j \in \llbracket N \rrbracket, i \neq j,
\end{align*}
by multiplying by the integer variables and their complements.
The system is then linearized by introducing auxiliary variables $y(u,v)$ for the product $uv$ of decision variables $u$ and $v$. For the product of the form $z_{r,s}^2$  we substitute $z_{r,s}$ after noting that $u^2=u$ for $u \in \{0,1\}$. For notational simplicity, we still write $y(z_{r,s},z_{r,s}) = z_{r,s}$.

This gives a lifted linear relaxation for \eqref{eqn:1D-FLP}, which is a mixed-integer analog to the first level of the Sherali-Adams construction;  denote this the \emph{lifted LP representation}. We may also consider adding a semidefinite constraint $M \succeq 0$ on the ``moment matrix'' $M$ given by
\[
    M_{I,J} = \begin{cases} 1 & I = J = \emptyset \\
                          z_{i,j} & I = (i,j), J = \emptyset \\
                          z_{i,j} & I = \emptyset, J = (i,j) \\
                          y(z_{i,j},z_{r,s}) & \text{o.w.} \: (I=(i,j),J=(r,s))
            \end{cases} \quad \forall I,J \in \{(r,s) \in \llbracket N \rrbracket^2 : r \neq s\} \cup \{\emptyset\}.
\]
We denote this the \emph{lifted SDP representation}, and note that it is reminiscent of the first level of the Lov\'{a}sz-Schrijver and Lasserre hierarchies, but in a mixed-integer setting.

We will show that applying either lifted representation to \eqref{eqn:1D-FLP} does not improve the dual bound beyond the closed-form combinatorial bound $\omega_2$.

\begin{proposition} \label{prop:level-one-hierarchy}
    The lifted LP and SDP representations of formulation \eqref{eqn:1D-FLP}, augmented with inequality \eqref{eqn:obj-cut}, have optimal cost equal to $\omega_2$.
\end{proposition}

\proof{}
    See Appendix~\ref{app:level-one-proof}.
\endproof

%

We note that the hierarchies presented in the previous section by no means preclude other, stronger extended formulations that may be constructed in a problem-specific manner. Indeed, such a specialized SDP formulation exists for the FLP, although we show in the following section that it does not improve on the closed form bound $\omega_2$.

\subsection{Ad-hoc SDP formulations for the FLP}
We now turn our attention to lifted representations of the FLP. \cite{Takouda:2005} present an SDP formulation for the entire FLP, which expresses the non-overlapping constraints (\ref{eqn:2D-FLP-6},\ref{eqn:2D-FLP-7},\ref{eqn:2D-FLP-11}) with complementarity conditions over $-1/1$ variables, whose integrality are then relaxed in the standard way (see e.g. Chapter 10 of \cite{Conforti:2014}). However, we show that the relaxation of this ad-hoc formulation produces a dual bound equal to the level-2 combinatorial bound $\omega_2$.

\begin{proposition} \label{prop:takouda-sdp}
    The relaxation of the SDP formulation from \cite{Takouda:2005} has optimal objective value equal to $\omega_2$.
\end{proposition}
\proof{}
    See Appendix~\ref{app:takouda-sdp}.
\endproof

\section{Computational results}\label{comp_sec}
Having compared the strength of our combinatorial bounding scheme against existing standard methods in previous sections, we will now perform a computational comparison. Specifically, we will be interested in comparing the computational effort (as measured in running time) needed to construct the same dual bound value. Specifically, we compare our combinatorial bounding scheme against: (i) The dual bound provided by a MIP solver, (ii) the lifted LP representation, and (iii) the lifted SDP representation.

First, we will compare the second level of the combinatorial bound against the lifted LP and SDP representations. Since we have shown in Proposition \ref{prop:level-one-hierarchy} that the combinatorial scheme is no worse, we will compare the running time of the two approaches to construct the same bound.

Secondly, we will perform a comparison of the combinatorial bound and the dual bound from a state-of-the-art MIP solver. Since there does not exist a nice relationship between the quality of the bounds for the two, we instead will construct the combinatorial bound for different levels, and then allot the same running time for the MIP solver. In this way, we will be comparing the quality of the bounds between the two approaches, given the same computational budget. We compare the relative gap for both approaches with respect to the best known feasible solution found by the MIP solver (optimal or not).

We perform the computational trials on a OS X machine with a 2.7GHz Intel Core i5-5257U processor and 8GB of RAM. We use Gurobi 6.0.4 as the IP solver and Mosek 7.1.0.31 as the SDP solver. We use the JuMP modeling language \cite{Dunning:2015a} in the Julia programming language \cite{Bezanson:2012} for algebraic modeling and scripting the computational trials.

Finally, we note that we make no attempt in these computational experiments of computing the partial dual bounds $\gamma(S(C))$ in parallel, which could potentially offer a speed-up over our serial implementation of the combinatorial dual bounding scheme. Per the note in Section \ref{comb_sec}, we compute in closed form all partial bounds $\gamma(C(S))$ for sets $|C|=2$; for all other subsets, we also use Gurobi and the same MIP formulation (\eqref{eqn:1D-FLP} or \eqref{eqn:2D-FLP}) to produce the partial bound value. For consistency of implementation, we still compute the bound $\omega_2$ via \eqref{eqn:1D-FLP-bounding-scheme} using Gurobi, even though it also is available in closed form as noted in Section \ref{ss:lp-relaxation}.

\subsection{1D-FLP computational results}
The benchmark instances are primarily culled from \cite{Amaral:2006,Amaral:2008}, with one large-scale (AMI33) instance taken as an adapted FLP instance from the MCNC instance set (see \cite{Huchette:2015}). Note that the number in the names indicate the number of components in the given instance.

\begin{table}[]
\centering
\begin{tabular}{r|@{\hskip 0.5em}r@{\hskip 0.5em}r@{\hskip 0.5em}r@{\hskip 0.5em}r@{\hskip 0.5em}r@{\hskip 0.5em}r@{\hskip 0.5em}r@{\hskip 0.5em}r@{\hskip 0.5em}r@{\hskip 0.5em}r@{\hskip 0.5em}r}
     &    S8 &  S8H &    S9 &   S9H &    S10 &    S11 &   LW11 &     P15 &     P17 &     P18 & AMI33 \\\hline
 CB  &  1.9e-5 & 1.9e-5 & 2.4e-5 & 2.3e-5 & 2.9e-5 & 3.6e-5 & 3.4e-5 & 6.0e-5 & 7.6e-5 & 8.5e-5 & 3.2e-4 \\
 LP  &    1.51 & 0.54 &     0.90 &   0.96 &   1.75 &   2.81 &   2.59 &  13.43 &  27.12 &  42.18 &    T/O \\
 SDP &   15.13 & 9.25 &    47.30 &  42.02 & 130.40 & 310.63 & 325.81 &    T/O &    T/O &    T/O &    T/O
\end{tabular}
\caption{Running time (in seconds) for the combinatorial bounding scheme (CB), LP representation, and SDP representation on the 1D-FLP benchmark instances. T/O indicates that the solver did not terminate in 10 minutes.}
\label{table:lifted}
\end{table}

In Table \ref{table:lifted} we compare the running time of the combinatorial bounding scheme for $k=2$ and the first level of the lifted LP and SDP representations. Recall from Proposition \ref{prop:level-one-hierarchy} that the combinatorial bound is as least as good as that from the lifted representations. Furthermore, we see that the combinatorial bound can be computed in less than $1e-4$ seconds for all instances, since the partial bounds are available in closed form. Contrastingly, the LP representation takes on the order of seconds for most of the instances, and the SDP representation roughly an order of magnitude longer. In particular, we see that the SDP representation is too large to be solved in under 10 minutes for 15 intervals or more, and that the LP representation also is too large for the largest instance in the test set.

We also note that we attempted to construct the second level of our lifted representation (in the sense of Lov\'{a}sz-Schrijver and Lasserre), but were unable to solve the smallest 1D-FLP instance (or even build many of the instances in memory). Our assessment is that the lifted representations, while theoretically interesting, do not provide an off-the-shelf practical approach for this problem as they are too large and not sufficiently strong.

Our second set of trials compares the combinatorial scheme against a MIP solver (Gurobi) given the same amount of time. The results are shown in the left half of Table \ref{table:FLP-budget} for the combinatorial scheme $\omega_k$ run for $k \in \{2,3,4,5\}$. We see that the combinatorial scheme does better on average, although there are two instances where Gurobi is able to find the optimal solution very quickly. We see that it is relatively cheap to compute the combinatorial bound for smaller instances or for smaller values of $k$, although for the larger instances the bound can get expensive for larger $k$. However, Gurobi struggles even more on these difficult instances, and hence we see noticeable gains from including higher levels of the scheme. We note that, since evaluating $\omega_2$ is is very inexpensive with the combinatorial bounding scheme, Gurobi is unable to return \emph{any} dual bound on the problem in the same amount of time.

\begin{table}
    \centering
    \begin{minipage}[b]{.45\linewidth} \centering
    {\scriptsize
        \begin{tabular}{cc|@{\hskip 1em}r@{\hskip 1em}r@{\hskip 1em}|r}
         & & CB & MIP & Time \\ \hline
        \multirow{4}{*}{S8} & 2 & \cellcolor{lightgray} 54.06 & -- & 1.9e-5 \\
         & 3 & 46.13 & \cellcolor{lightgray} 43.57 & 0.17 \\
         & 4 & 37.27 & \cellcolor{lightgray}  6.37 & 0.53 \\
         & 5 & 29.01 & \cellcolor{lightgray}  0.00 & 0.98 \\ \hline

        \multirow{4}{*}{S8H} & 2 & \cellcolor{lightgray} 63.67 & -- & 1.9e-5 \\
         & 3 & \cellcolor{lightgray} 53.17 & 72.64 & 0.17 \\
         & 4 & \cellcolor{lightgray} 42.80 & 43.49 & 0.53 \\
         & 5 & \cellcolor{lightgray} 32.25 & 41.77 & 0.97 \\ \hline

        \multirow{4}{*}{S9} & 2 & \cellcolor{lightgray} 57.99 & -- & 2.4e-5 \\
         & 3 & \cellcolor{lightgray} 49.49 & 55.08 & 0.31 \\
         & 4 & 41.28 & \cellcolor{lightgray} 20.46 & 0.86 \\
         & 5 & 32.80 & \cellcolor{lightgray}  0.00 & 2.01 \\ \hline

        \multirow{4}{*}{S9H} & 2 & \cellcolor{lightgray} 67.64 & -- & 2.3e-5 \\
         & 3 & \cellcolor{lightgray} 58.20 & 71.53 & 0.25 \\
         & 4 & \cellcolor{lightgray} 49.11 & 53.69 & 0.81 \\
         & 5 & \cellcolor{lightgray} 39.61 & 52.91 & 1.96 \\ \hline

        \multirow{4}{*}{S10} & 2 & \cellcolor{lightgray} 58.71 & -- & 2.9e-5 \\
         & 3 & \cellcolor{lightgray} 50.39 & 54.25 & 0.40 \\
         & 4 & 43.46 & \cellcolor{lightgray} 37.40 & 1.35 \\
         & 5 & \cellcolor{lightgray} 35.54 & 35.80 & 3.45 \\ \hline

        \multirow{4}{*}{S11} & 2 & \cellcolor{lightgray} 63.30 & -- & 3.6e-5 \\
         & 3 & \cellcolor{lightgray} 56.38 & 68.43 & 0.47 \\
         & 4 & 49.58 & \cellcolor{lightgray} 45.03 & 2.05 \\
         & 5 & 42.82 & \cellcolor{lightgray} 39.60 & 6.28 \\ \hline

        \multirow{4}{*}{LW11} & 2 & \cellcolor{lightgray} 63.30 & -- & 3.4e-5 \\
         & 3 & \cellcolor{lightgray} 56.38 & 65.55 & 0.55 \\
         & 4 & 49.58 & \cellcolor{lightgray} 48.12 & 2.03 \\
         & 5 & \cellcolor{lightgray} 42.82 & 43.54 & 6.24 \\ \hline

        \multirow{4}{*}{P15} & 2 & \cellcolor{lightgray} 69.18 & -- & 6.0e-5 \\
         & 3 & \cellcolor{lightgray} 64.01 & 77.34 & 1.11 \\
         & 4 & \cellcolor{lightgray} 59.06 & 71.32 & 6.91 \\
         & 5 & \cellcolor{lightgray} 53.95 & 55.84 & 32.37 \\ \hline

        \multirow{4}{*}{P17} & 2 & \cellcolor{lightgray} 72.61 & -- & 7.6e-5 \\
         & 3 & \cellcolor{lightgray} 68.14 & 82.17 & 1.78 \\
         & 4 & \cellcolor{lightgray} 63.32 & 72.66 & 12.76 \\
         & 5 & \cellcolor{lightgray} 58.34 & 65.17 & 63.03 \\ \hline

        \multirow{4}{*}{P18} & 2 & \cellcolor{lightgray} 74.27 & -- & 8.5e-5 \\
         & 3 & \cellcolor{lightgray} 69.91 & 83.34 & 2.15 \\
         & 4 & \cellcolor{lightgray} 65.33 & 74.53 & 15.49 \\
         & 5 & \cellcolor{lightgray} 60.68 & 67.36 & 85.28 \\ \hline

        \multirow{4}{*}{AMI33} & 2 & \cellcolor{lightgray} 83.45 & -- & 3.2e-4 \\
         & 3 & \cellcolor{lightgray} 80.42 & 95.86 &  17.24 \\
         & 4 & \cellcolor{lightgray} 77.28 & 86.31 &  207.68 \\
         & 5 & \cellcolor{lightgray} 74.64 & 85.15 & 2262.69
    \end{tabular}
    }
    \end{minipage} \hspace{2em}
    \begin{minipage}[b]{.45\linewidth} \centering
    {\scriptsize
        \begin{tabular}{cc|@{\hskip 1em}r@{\hskip 1em}r@{\hskip 1em}|r}
         & & CB & MIP & Time \\ \hline

        \multirow{4}{*}{A9} & 2 & \cellcolor{lightgray} 58.36 & -- & 3.1e-5 \\
         & 3 & \cellcolor{lightgray} 50.03 & 79.55 & 0.88 \\
         & 4 & \cellcolor{lightgray} 42.86 & 76.48 & 6.72 \\
         & 5 & \cellcolor{lightgray} 34.24 & 67.61 & 37.59 \\ \hline

        \multirow{4}{*}{B9} & 2 & \cellcolor{lightgray} 61.60 & -- & 3.0e-5 \\
        & 3 & \cellcolor{lightgray} 57.97 & 74.18 & 0.46 \\
        & 4 & \cellcolor{lightgray} 48.02 & 65.79 & 2.71 \\
        & 5 & \cellcolor{lightgray} 42.21 & 54.92 & 12.27 \\ \hline

        \multirow{4}{*}{X10} & 2 & \cellcolor{lightgray} 56.16 & -- &  3.8e-5 \\
        & 3 & \cellcolor{lightgray} 49.06 & 71.70 &  1.44 \\
        & 4 & \cellcolor{lightgray} 41.43 & 50.22 & 10.37 \\
        & 5 & 34.90 & \cellcolor{lightgray} 31.72 & 70.14 \\ \hline

        \multirow{4}{*}{C10} & 2 & \cellcolor{lightgray} 44.03 & -- & 3.4e-5 \\
        & 3 & \cellcolor{lightgray} 40.19 & 99.61 & 0.23 \\
        & 4 & \cellcolor{lightgray} 35.79 & 80.86 & 1.90 \\
        & 5 & \cellcolor{lightgray} 33.54 & 67.10 & 16.72 \\ \hline

        \multirow{4}{*}{H11} & 2 & \cellcolor{lightgray} 51.51 & -- & 3.8e-5 \\
        & 3 & \cellcolor{lightgray} 43.22 & 73.96 & 1.50 \\
        & 4 & \cellcolor{lightgray} 36.01 & 56.63 & 13.81 \\
        & 5 & \cellcolor{lightgray} 30.14 & 37.18 & 112.11 \\ \hline

        \multirow{4}{*}{B12} & 2 & \cellcolor{lightgray} 55.43 & -- & 4.7e-5 \\
        & 3 & \cellcolor{lightgray} 47.93 & 68.18 & 0.39 \\
        & 4 & \cellcolor{lightgray} 43.52 & 52.57 & 3.11 \\
        & 5 & \cellcolor{lightgray} 37.07 & 46.57 & 19.71 \\ \hline

        \multirow{4}{*}{B13} & 2 & \cellcolor{lightgray} 63.21 & -- & 5.7e-5 \\
        & 3 & \cellcolor{lightgray} 52.62 & 75.62 & 2.45 \\
        & 4 & \cellcolor{lightgray} 44.57 & 75.54 & 23.79 \\
        & 5 & \cellcolor{lightgray} 38.37 & 61.54 & 203.24 \\ \hline

        \multirow{4}{*}{B14} & 2 & \cellcolor{lightgray} 73.29 & -- & 6.1e-5 \\
        & 3 & \cellcolor{lightgray} 67.55 & 84.10 & 2.56 \\
        & 4 & \cellcolor{lightgray} 62.48 & 80.03 & 27.85 \\
        & 5 & \cellcolor{lightgray} 58.44 & 71.74 & 276.85 \\ \hline

        \multirow{4}{*}{B15} & 2 & \cellcolor{lightgray} 43.77 & -- & 6.9e-5 \\
        & 3 & \cellcolor{lightgray} 41.80 & 98.49 & 0.65 \\
        & 4 & \cellcolor{lightgray} 34.97 & 95.73 & 8.49 \\
        & 5 & \cellcolor{lightgray} 32.50 & 88.26 & 149.10 \\ \hline

        \multirow{4}{*}{A20} & 2 & \cellcolor{lightgray} 77.74 & -- & 1.2e-4 \\
        & 3 & \cellcolor{lightgray} 75.50 & 87.55 & 2.99 \\
        & 4 & \cellcolor{lightgray} 72.55 & 86.47 & 52.10 \\
        & 5 & \cellcolor{lightgray} 70.02 & 85.27 & 1069.01 \\ \hline

        \multirow{4}{*}{A20M} & 2 & \cellcolor{lightgray} 78.03 & -- & 1.2e-4 \\
        & 3 & \cellcolor{lightgray} 75.61 & 90.16 & 3.27 \\
        & 4 & \cellcolor{lightgray} 72.45 & 90.16 & 57.48 \\
        & 5 & \cellcolor{lightgray} 70.75 & 90.16 & 1225.73
    \end{tabular}
    }
    \end{minipage}

    \caption{Comparison of the combinatorial bounding scheme (CB) and a MIP solver (MIP) on the 1D-FLP instances (Left) and FLP instances (Right). Note that the numbers in the benchmark names indicate the number of components/intervals. We run the CB scheme for $k \in \{2,3,4,5\}$ as denoted in the first column. We report the relative gap percentage ($100(UB-LB)/UB$), and the time budget allotted to each method (in seconds). Finally, we color in gray the method that produces the better (i.e. smaller) relative gap.}
    \label{table:FLP-budget}
\end{table}

However, we stress that these results are restricted to formulation \eqref{eqn:1D-FLP} for the 1D-FLP, as they do not compare against the bounds from \cite{Amaral:2009} and \cite{Anjos:2005,Anjos:2008}, which offer promising alternative dual bounds for the 1D-FLP. 

\subsection{FLP computational results}

\begin{table}[]
\centering
\begin{tabular}{r|@{\hskip 0.5em}r@{\hskip 0.5em}r@{\hskip 0.5em}r@{\hskip 0.5em}r@{\hskip 0.5em}r@{\hskip 0.5em}r@{\hskip 0.5em}r@{\hskip 0.5em}r@{\hskip 0.5em}r@{\hskip 0.5em}r@{\hskip 0.5em}r}
     &    A9 &       B9 &    X10 &    C10 &    H11 &    B12 &    B13 &    B14 &    B15 &    A20 &   A20M \\\hline
 CB  &  3.1e-5 & 3.0e-5 & 3.8e-5 & 3.4e-5 & 3.8e-5 & 4.7e-5 & 5.7e-5 & 6.1e-5 & 6.9e-5 & 1.2e-4 & 1.2e-4 \\
 SDP &   83.23 &  98.95 & 310.57 & 327.44 &    T/O &    T/O &    T/O &    T/O &    T/O &    T/O &    T/O
\end{tabular}
\caption{Running time (in seconds) for the combinatorial bounding scheme (CB) and the SDP formulation from \cite{Takouda:2005} (SDP) on the FLP benchmark instances. T/O indicates that the solver did not terminate in 10 minutes.}
\label{table:FLP-SDP-results}
\end{table}

Our benchmark instance collection for the FLP is the same 11 instances from the literature used in \cite{Huchette:2015}. In Table \ref{table:FLP-SDP-results}, we present the running times for computing $\omega_2$ from the combinatorial scheme and for the lifted SDP formulation in \cite{Takouda:2005}; again, in Proposition \ref{prop:takouda-sdp} we have shown that the two approaches produce the same bound value. We do not include the lifted LP or lifted SDP formulations in the comparison, since \emph{none} of the instances were able to build in memory in under 10 minutes, much less solve to optimality. We see that the SDP formulation from \cite{Takouda:2005} fares slightly better, but is still only able to solve to optimality on the smaller instances ($\leq 10$ components), and in orders of magnitude more time than the combinatorial scheme requires.

In the right half of Table \ref{table:FLP-budget} we record the computational comparisons between the combinatorial bounding scheme and the dual bound produced by Gurobi. We see that the combinatorial bounding scheme uniformly returns a superior dual bound (with a single exception), sometimes by a significant margin. In particular, we highlight the last instance, A20M, a particularly difficult 20 component instance. The MIP solver is unable to make any noticeable progress on the dual bound of this problem, even after almost 10 minutes (this trend continues if the solver is left running for longer periods as well). In contrast, the combinatorial bounding scheme is able to make noticeable improvements to the gap at successively higher levels for $k$.

We are also cautiously optimistic that the gaps for the difficult instances would be significantly better (i.e. smaller) if the relative gap were computed with the true (unknown) optimal cost. This is because there is a wide body of literature constructing high-performing heuristics for the FLP that are likely to be of higher quality than those produced by the IP solver we used for this comparison.


\section*{Acknowledgements}
This material is based upon work supported by the National Science Foundation Graduate Research Fellowship under Grant No. 1122374 and Grant CMMI-1351619.

\bibliographystyle{ormsv080}
{\footnotesize \bibliography{master}}

\newpage
\begin{appendix}

\section{Proof of Proposition \ref{prop:level-one-hierarchy}} \label{app:level-one-proof}
\proof{}
    See that \eqref{eqn:obj-cut} immediately gives that the optimal cost can be no less than $\omega_2$; we just need to show that it can also be no greater. The linearized system (sans semidefinite constraint) is, after some rearranging, given by
    \begin{subequations} \label{eqn:sherali-adams-1}
    \begin{align}
        d_{i,j} - c_i + c_j \geq y(d_{i,j},z_{r,s}) - y(c_i,z_{r,s}) + y(c_j,z_{r,s}) \geq 0 \label{eqn:sherali-adams-1-1} \\
        d_{i,j} + c_i - c_j \geq y(d_{i,j},z_{r,s}) + y(c_i,z_{r,s}) - y(c_j,z_{r,s}) \geq 0 \label{eqn:sherali-adams-1-2} \\
        L(1-z_{i,j}-z_{r,s}+y(z_{i,j},z_{r,s})) - c_i + c_j - (\ell_i+\ell_j)(1-z_{r,s}) \geq - y(c_i,z_{r,s}) + y(c_j,z_{r,s}) \label{eqn:sherali-adams-1-3} \\
        Lz_{r,s} - Ly(z_{i,j},z_{r,s}) - y(c_i,z_{r,s}) + y(c_j,z_{r,s}) - (\ell_i+\ell_j)z_{r,s} \geq 0 \label{eqn:sherali-adams-1-4} \\
        y(z_{i,j},z_{r,s}) + y(z_{j,i},z_{r,s}) = z_{r,s} \label{eqn:sherali-adams-1-5} \\
        z_{i,j} + z_{j,i} = 1 \label{eqn:sherali-adams-1-6} \\
        c_i \geq y(c_i,z_{r,s}) \geq 0 \label{eqn:sherali-adams-1-7} \\
        L - c_i \geq Lz_{r,s} - y(c_i,z_{r,s}) \geq 0 \label{eqn:sherali-adams-1-8} \\
        d_{i,j} - (\ell_i+\ell_j) \geq y(d_{i,j},z_{r,s}) - (\ell_i+\ell_j)z_{r,s} \geq 0 \label{eqn:sherali-adams-1-9} \\
        z_{i,j} \geq y(z_{i,j},z_{r,s}) \geq 0 \label{eqn:sherali-adams-1-10} \\
        1 - z_{i,j} \geq z_{r,s} - y(z_{i,j},z_{r,s}) \geq 0, \label{eqn:sherali-adams-1-11}
    \end{align}
    \end{subequations}
    for all $(i,j),(r,s) \in P(\llbracket N \rrbracket)$. The proposed feasible solution is
    \begin{align*}
        \hat{c}_i &\leftarrow L/2 \\
        \hat{d}_{i,j} &\leftarrow \ell_i+\ell_j \\
        \hat{z}_{i,j} &\leftarrow 1/2 \\
        \hat{y}(c_i,z_{r,s}) &\leftarrow \begin{cases}
                \frac{1}{2}(L/2 - \ell_s) & r = i \\
                \frac{1}{2}(L/2 + \ell_r) & s = i \\
                L/4          & \text{o.w.}
            \end{cases} \\
        \hat{y}(d_{i,j},z_{r,s}) &\leftarrow \frac{1}{2}(\ell_i+\ell_j) \\
        \hat{y}(z_{i,j},z_{r,s}) &\leftarrow \begin{cases}
                1/2 & i = r, j = s \\
                0   & i = s, j = r \\
                1/4 & \text{o.w.}
            \end{cases}
    \end{align*}
    for all $(i,j), (r,s) \in P(\llbracket N \rrbracket)$. It is clear that this has objective value equal to $\omega_2$, so it suffices to show feasibility w.r.t. \eqref{eqn:sherali-adams-1}, which may be verified. In particular, we show \eqref{eqn:sherali-adams-1-1}, \eqref{eqn:sherali-adams-1-3}, \eqref{eqn:sherali-adams-1-4}, \eqref{eqn:sherali-adams-1-5}, \eqref{eqn:sherali-adams-1-7}, and \eqref{eqn:sherali-adams-1-8}, as the rest are immediate or follow analogously.

    \paragraph{\eqref{eqn:sherali-adams-1-1}}
        First, we observe that
        \[
            -\frac{1}{2}(\ell^s_i+\ell^s_j) \leq -\hat{y}(c_i,z_{r,s}) + \hat{y}(c_j,z_{r,s}) \leq \frac{1}{2}(\ell^s_i+\ell^s_j),
        \]
        which implies the desired result
        \[
            \hat{d}_{i,j} - \hat{c}_i + \hat{c}_j = \ell^s_i + \ell^s_j \geq \hat{y}(d_{i,j},z_{r,s}) - \hat{y}(c_i,z_{r,s}) + \hat{y}(c_j,z_{r,s}) \geq 0.
        \]
    \paragraph{\eqref{eqn:sherali-adams-1-3}}
        First, we see that the left-hand side reduces to
        \[
            L\hat{y}(z_{i,j},z_{r,s}) - \frac{1}{2}(\ell_i+\ell_j).
        \]
        In the case that $(i,j) = (r,s)$, we need that
        \[
            L/2 - \frac{1}{2}(\ell_i+\ell_j) \geq -(L/4 - \ell_j/2) + (L/4 + \ell_i/2),
        \]
        which holds true.
        In the case that $(i,j) = (s,r)$, we need that
        \[
            -\frac{1}{2}(\ell_i+\ell_j) \geq -(L/4 + \ell_j/2) + (L/4 - \ell_i/2),
        \]
        which also holds true. In the case that $i=r$ and $j\neq s$,
        \[
            L/4  - \frac{1}{2}(\ell_i+\ell_j) \geq -(L/4-\ell_s/2) + L/4,
        \]
        which holds from feasibility. For $i=s$, $j \neq r$, we get
        \[
            L/4  - \frac{1}{2}(\ell_i+\ell_j) \geq -(L/4+\ell_r/2) + L/4;
        \]
        for $i\neq r$, $j=s$, we get
        \[
            L/4  - \frac{1}{2}(\ell_i+\ell_j) \geq -L/4 + (L/4+\ell_r/2);
        \]
        for $i\neq s$, $j=r$, we get
        \[
            L/4  - \frac{1}{2}(\ell_i+\ell_j) \geq -L/4 + (L/4-\ell_s/2);
        \]
        for $i,j \notin \{r,s\}$, we get
        \[
            L/4  - \frac{1}{2}(\ell_i+\ell_j) \geq -L/4 + L/4,
        \]
        all true statements. This exhausts all possible cases.

    \paragraph{\eqref{eqn:sherali-adams-1-4}}
        Similarly as to (\ref{eqn:sherali-adams-1}c), this reduces to showing that
        \[
            L/2 - \frac{1}{2}(\ell_i+\ell_j) - \hat{y}(c_i,z_{r,s}) + \hat{y}(c_j,z_{r,s}) \geq 0.
        \]
        In the case that $(i,j) = (r,s)$, we need that
        \[
            L/2 - \frac{1}{2}(\ell_i+\ell_j) -(L/4 - \ell_j/2) + (L/4 + \ell_i/2) \geq 0,
        \]
        which holds true.
        In the case that $(i,j) = (s,r)$, we need that
        \[
            L/2 - \frac{1}{2}(\ell_i+\ell_j) -(L/4 + \ell_j/2) + (L/4 - \ell_i/2) \geq 0,
        \]
        which also holds true. In the case that $i=r$ and $j\neq s$,
        \[
            L/2  - \frac{1}{2}(\ell_i+\ell_j) - (L/4-\ell_s/2) + L/4 \geq 0,
        \]
        which holds from feasibility. For $i=s$, $j \neq r$, we get
        \[
            L/2  - \frac{1}{2}(\ell_i+\ell_j) - (L/4+\ell_r/2) + L/4 \geq 0;
        \]
        for $i\neq r$, $j=s$, we get
        \[
            L/2  - \frac{1}{2}(\ell_i+\ell_j) - L/4 + (L/4+\ell_r/2) \geq 0;
        \]
        for $i\neq s$, $j=r$, we get
        \[
            L/2  - \frac{1}{2}(\ell_i+\ell_j) - L/4 + (L/4-\ell_s/2) \geq 0;
        \]
        for $i,j \notin \{r,s\}$, we get
        \[
            L/2  - \frac{1}{2}(\ell_i+\ell_j) - L/4 + L/4 \geq 0,
        \]
        all true statements. This exhausts all possible cases.

    \paragraph{\eqref{eqn:sherali-adams-1-5}}
        If $(i,j) = (s,r)$, then
        \[
            \hat{y}(z_{i,j},z_{j,i}) + \hat{y}(z_{j,i},z_{j,i}) = 0 + 1/2 = 1/2 = \hat{z}_{r,s}.
        \]
        Similarly if $(i,j) = (r,s)$. Otherwise, we have
        \[
            \hat{y}(z_{i,j},z_{r,s}) + \hat{y}(z_{j,i},z_{r,s}) = 1/4 + 1/4 = 1/2 = \hat{z}_{r,s}.
        \]

    \paragraph{\eqref{eqn:sherali-adams-1-7}}
        Since $2\ell_i \leq L$ from feasibility,
        \[
            0 \leq \hat{y}(c_i,z_{r,s}) \leq L/2,
            \]
        which gives the result.

    \paragraph{\eqref{eqn:sherali-adams-1-8}}
        Using the fact from (\ref{eqn:sherali-adams-1}g), we get that
        \begin{gather*}
            L - \hat{c}_i = L/2 \\
            0 \leq L\hat{z}_{r,s} - \hat{y}(c_i,z_{r,s}) \leq L/2,
        \end{gather*}
        which gives the result.

    To show that the moment matrix satisfies the semidefinite constraint $M \succeq 0$, consider that for the proposed solution,
    \[
        M_{I,J} = \begin{cases}
            1   & I = J = \emptyset \\
            1/2 & I = J \neq \emptyset \\
            1/2 & I = \emptyset \text{ or } J = \emptyset \\
            0   & I = (i,j), J = (j,i) \\
            1/4 & \text{o.w.}
        \end{cases} \quad \forall I,J \in \{(r,s) \in \llbracket N \rrbracket^2 : r \neq s\} \cup \{\emptyset\}.
    \]
    We will derive a congruence transformation between this matrix and a diagonal matrix via elementary row and column operations; then, showing that the diagonal matrix is positive semidefinite gives the result.

    Produce $M'$ by adding $-1/2$ times the first row ($\emptyset$) to each row $I \neq \emptyset$. Then
    \[
        M'_{I,J} = \begin{cases}
            1    & I = J = \emptyset \\
            1/4  & I = J \neq \emptyset \\
            1/2  & I = \emptyset, J \neq \emptyset \\
            -1/4 & I = (i,j), J = (j,i) \\
            0    & \text{o.w.}
        \end{cases} \quad \forall I,J \in \{(r,s) \in \llbracket N \rrbracket^2 : r \neq s\} \cup \{\emptyset\}.
    \]
    Now add row $(i,j)$ to row $(j,i)$ (where $i < j$) to produce the matrix
    \[
        M''_{I,J} = \begin{cases}
            1    & I = J = \emptyset \\
            1/4  & I = J = (i,j), i < j \\
            1/2  & I = \emptyset, J \neq \emptyset \\
            -1/4 & I = (i,j), J = (j,i), i < j \\
            0    & \text{o.w.}
        \end{cases} \quad \forall I,J \in \{(r,s) \in \llbracket N \rrbracket^2 : r \neq s\} \cup \{\emptyset\}.
    \]
    We now perform the same series of elementary operations, but to the columns instead of rows. By adding $1/2$ times the first column ($\emptyset$) to each column $J \neq \emptyset$, we get
    \[
        M'''_{I,J} = \begin{cases}
            1    & I = J = \emptyset \\
            1/4  & I = J = (i,j), i < j \\
            -1/4 & I = (i,j), J = (j,i), i < j \\
            0    & \text{o.w.}
        \end{cases} \quad \forall I,J \in \{(r,s) \in \llbracket N \rrbracket^2 : r \neq s\} \cup \{\emptyset\}.
    \]
    Finally, add column $(i,j)$ to column $(j,i)$ (where $i < j)$ to produce the diagonal matrix
    \[
        M''''_{I,J} = \begin{cases}
            1    & I = J = \emptyset \\
            1/4  & I = J = (i,j), i < j \\
            0    & \text{o.w.}
        \end{cases} \quad \forall I,J \in \{(r,s) \in \llbracket N \rrbracket^2 : r \neq s\} \cup \{\emptyset\}.
    \]

    This produces a diagonal matrix (under an appropriate ordering of the rows and columns) whose diagonal entries are non-negative. Therefore, $M''''$ is positive semidefinite, and therefore by Sylvester's law of inertia \cite{Horn:1990}, $M$ is as well, since we have produced $M''''$ via a congruence transform $M'''' = HMH^T$, where $H$ is the invertible matrix describing the elementary row operations. This yields feasibility w.r.t. the lifted SDP representation.
\endproof

\section{Proof of Proposition \ref{prop:takouda-sdp}} \label{app:takouda-sdp}
\proof{}
The notation for the FLP in \cite{Takouda:2005} is different than in this work; we keep the original notation for simplicity. In particular, the floor is shifted from $[0,L^x]\times[0,L^y]$ to $[-L^x/2,L^x/2] \times [-L^y/2,L^y/2]$, the desired area is notated with $a_i$ instead of $\alpha_i$, and $\ell_i$ now denotes the widths of component $i$, as opposed the half-width (the bounds $lb_i$ and $ub_i$ are now for the width as well, accordingly). We note that we omit the aspect ratio constraints included in \cite{Takouda:2005} in lieu of enforcing them via bounds on $\ell$ (see \cite{Castillo:2005} or \cite{Huchette:2015} for details). Additionally, we omit the $p-q$ symmetry-breaking constraints mentioned in Section 4.3 of \cite{Takouda:2005}, as they are identical to the symmetry-breaking constraints we have omitted from our other formulations. We now present their proposed model as:
\begin{subequations} \label{eqn:takouda}
\begin{gather}
\min \sum_{(i,j)\in P(\llbracket N \rrbracket)} p_{i,j}(d_{i,j}^x+d_{i,j}^y) \\
\operatorname{s.t.} \begin{pmatrix} \ell^x_i & \sqrt{a_i} \\ \sqrt{a_i} & \ell^y_i \end{pmatrix} \succeq 0 \label{eqn:takouda-1} \\
\begin{pmatrix} \frac{1}{2}(L^s-\ell^s_i) & c^s_i \\ c^s_i & \frac{1}{2}(L^s-\ell^s_i) \end{pmatrix} \succeq 0 \label{eqn:takouda-2} \\
d_{i,j}^x \geq \frac{1}{2}(\ell^x_i+\ell^x_j) - \frac{L^x}{2}(1-\sigma_{i,j}) \label{eqn:takouda-3} \\
d_{i,j}^y \geq \frac{1}{2}(\ell^y_i+\ell^y_j) - \frac{L^y}{2}(1+\sigma_{i,j}) \label{eqn:takouda-4} \\
d_{i,j}^s - 2S_{i,j}^s = c^s_j - c^s_i \label{eqn:takouda-5} \\
S_{i,j}^s \geq 0 \label{eqn:takouda-6} \\
S_{i,j}^s + c^s_j - c^s_i \geq 0 \label{eqn:takouda-7} \\
S_{i,j}^x \leq \frac{1}{4}\left(3 - \sigma_{i,j} - \alpha_{i,j} - \sigma_{i,j}\alpha_{i,j}\right)\left(L^x-\frac{1}{2}(lb^x_i+lb^x_j)\right) \label{eqn:takouda-8} \\
S_{i,j}^x + c^x_j - c^x_i \leq \frac{1}{4}\left( 3 - \sigma_{i,j} + \alpha_{i,j} + \sigma_{i,j}\alpha_{i,j} \right)\left(L^x-\frac{1}{2}(lb^x_i+lb^x_j)\right) \label{eqn:takouda-9} \\
S_{i,j}^y \leq \frac{1}{4}\left( 3 + \sigma_{i,j} - \alpha_{i,j} + \sigma_{i,j}\alpha_{i,j} \right)\left(L^y-\frac{1}{2}(lb^y_i+lb^y_j)\right) \label{eqn:takouda-10} \\
S_{i,j}^y + c^y_j - c^y_i \leq \frac{1}{4}\left( 3 + \sigma_{i,j} + \alpha_{i,j} - \sigma_{i,j}\alpha_{i,j} \right)\left(L^y-\frac{1}{2}(lb^y_i+lb^y_j)\right) \label{eqn:takouda-11} \\
(\sigma_{i,j}+\sigma_{j,k})(\sigma_{i,j}+\sigma_{i,k})(\alpha_{i,j}+\alpha_{j,k})(\alpha_{i,j}-\alpha_{i,k}) = 0 \label{eqn:takouda-12} \\
d_{i,j}^x \geq \frac{1}{2}(lb^x_i+lb^x_j)(1+\sigma_{i,j}) \label{eqn:takouda-13} \\
d_{i,j}^y \geq \frac{1}{2}(lb^y_i+lb^y_j)(1-\sigma_{i,j}) \label{eqn:takouda-14} \\
c^x_i - c^x_j +\frac{1}{2}(\ell^x_i+\ell^x_j) \leq \frac{1}{2}(3-\sigma_{i,j}-\alpha_{i,j}-\sigma_{i,j}\alpha_{i,j})L^x \label{eqn:takouda-15} \\
c^x_j - c^x_i +\frac{1}{2}(l^x_i+l^x_j) \leq \frac{1}{2}(3-\sigma_{i,j}+\alpha_{i,j}+\sigma_{i,j}\alpha_{i,j})L^x \label{eqn:takouda-16} \\
c^y_i - c^y_j +\frac{1}{2}(l^y_i+l^y_j) \leq \frac{1}{2}(3+\sigma_{i,j}-\alpha_{i,j}+\sigma_{i,j}\alpha_{i,j})L^y \label{eqn:takouda-17} \\
c^y_j - c^y_i +\frac{1}{2}(l^y_i+l^y_j) \leq \frac{1}{2}(3+\sigma_{i,j}+\alpha_{i,j}-\sigma_{i,j}\alpha_{i,j})L^y \label{eqn:takouda-18} \\
lb^s_i \leq \ell^s_i \leq ub^s_i \label{eqn:takouda-19} \\
-1 \leq \sigma_{i,j}, \alpha_{i,j} \leq 1 \label{eqn:takouda-20} \\
\sigma_{i,j}, \alpha_{i,j} \in \{-1,1\} \label{eqn:takouda-21}
\end{gather}
\end{subequations}
for all distinct $i,j,k \in \llbracket N \rrbracket$ and $s \in \{x,y\}$. In words, the variables $c^s_i,\ell^s_i$, and $d^s_{i,j}$ serve the same functional role as in formulation \eqref{eqn:2D-FLP}. Additionally, there are auxiliary variables $S^s_{i,j}$ which correspond to the slack on the objective variable $d^s_{i,j}$ in the scenario that $\ell^s_j < \ell^s_i$. The discrete decision variables are
\[
    \sigma_{i,j} = \begin{cases}
        +1 & i \text{ and } j \text{ are separated along direction } x \\
        -1 & i \text{ and } j \text{ are separated along direction } y
    \end{cases}
\]
and
\[
    \alpha_{i,j} = \begin{cases}
        +1 & i \text{ precedes } j \text{ in the chosen direction} \\
        -1 & j \text{ precedes } i \text{ in the chosen direction}.
    \end{cases}
\]
The constraints \eqref{eqn:takouda-1} impose the area constraints, and \eqref{eqn:takouda-2} constrains the components to lie completely on the floor. Taken together, (\ref{eqn:takouda-3}-\ref{eqn:takouda-12}) simultaneously impose the non-overlap constraints and that the objective variables $d$ take the correct value. Constraints (\ref{eqn:takouda-13}-\ref{eqn:takouda-14}) are the B2 valid inequalities from \cite{Meller:1999}, which are akin to \eqref{eqn:obj-cut-2D-FLP}, and (\ref{eqn:takouda-15}-\ref{eqn:takouda-18}) are so-called $S$ valid inequalities derived in \cite{Castillo:2005}.

Note that this formulation still includes the complementarity constraints \eqref{eqn:takouda-12}; we will construct a direct solution to the relaxation of this formulation as-is, which will imply a rank-one feasible solution for the relaxation of this system, which in turn will imply a feasible solution for the final SDP formulation. Set
\begin{alignat*}{2}
\hat{c}^s_i &\leftarrow 0 &\quad \forall s \in \{x,y\}, i \in \llbracket N \rrbracket \\
\hat{\ell}^s_i &\leftarrow \max\{\sqrt{a_i},lb^s_i\} &\quad \forall s \in \{x,y\}, i \in \llbracket N \rrbracket \\
\hat{\sigma}_{i,j} &\leftarrow 0 &\quad \forall (i,j) \in P(\llbracket N \rrbracket) \\
\hat{\alpha}_{i,j} &\leftarrow 0 &\quad \forall (i,j) \in P(\llbracket N \rrbracket) \\
\hat{d}_{i,j}^x &\leftarrow \begin{cases} \frac{1}{2}(lb^x_i+lb^x_j) & lb^x_i+lb^x_j \leq lb^y_i+lb^y_j \\ 0 & o.w. \end{cases} &\quad \forall (i,j) \in P(\llbracket N \rrbracket) \\
\hat{d}_{i,j}^y &\leftarrow \begin{cases} \frac{1}{2}(lb^y_i+lb^y_j) & lb^x_i+lb^x_j > lb^y_i+lb^y_j \\ 0 & o.w. \end{cases} &\quad \forall (i,j) \in P(\llbracket N \rrbracket) \\
\hat{S}_{i,j}^s &\leftarrow \frac{1}{2}\hat{d}_{i,j}^s &\quad \forall s \in \{x,y\}, (i,j) \in P(\llbracket N \rrbracket).
\end{alignat*}

It can be verified that this solution satisfies the relaxation (\ref{eqn:takouda-1}-\ref{eqn:takouda-20}). In particular, we consider constraints \eqref{eqn:takouda-3}, \eqref{eqn:takouda-8}, and \eqref{eqn:takouda-15}; the rest are immediate, or use identical arguments.

\paragraph{\eqref{eqn:takouda-3}}
This follows from verifying that
\[
    \hat{d}^x_{i,j} \geq 0 \geq \frac{1}{2}(\ell^x_i+\ell^x_j) - \frac{L^x}{2}
\]
from assumption on the size of the floor $L^x$.

\paragraph{\eqref{eqn:takouda-8}}
This follows from showing the second inequality in
\[
    S^x_{i,j} \leq \frac{1}{4}(lb^x_i+lb^x_j) \leq \frac{3}{4}\left( L^x - \frac{1}{2}(lb^x_i+lb^x_j) \right),
\]
which follows from assumption on the floor size $L^x$.

\paragraph{\eqref{eqn:takouda-15}}
This follows from showing that
\[
    \frac{1}{2}(\sqrt{\alpha_i} + \sqrt{\alpha_j}) \leq \frac{1}{2}(ub^s_i+ub^s_j) \leq \frac{3}{2}L^x,
\]
which is satisfied given our assumptions on the size of the floor.

Furthermore, the formulation from \cite{Takouda:2005} linearizes the products of $\sigma$ and $\alpha$ variables that appear in \eqref{eqn:takouda}, and adds a semidefinite constraint of the form $vv^T \succeq 0$ on these products, where
\[
    v \defeq (1, \sigma_{1,2}, \ldots, \sigma_{N-1,N}, \alpha_{1,2},\ldots,\alpha_{N-1,N}, \sigma_{1,2}\sigma_{1,3},\ldots,\sigma_{N-2,N}\sigma_{N-1,N},\alpha_{1,2}\alpha_{1,3},\ldots,\alpha_{N-2,N}\alpha_{N-1,N}).
\]
However, we have feasibility with respect to this constraint directly from construction.


To see that the relaxation value cannot be less than the level-2 bound, sum \eqref{eqn:takouda-13} + \eqref{eqn:takouda-14} to get
\begin{align*}
    d_{i,j}^x + d_{i,j}^y &\geq \frac{1}{2}(lb^x_i + lb^x_j + lb^y_i + lb^y_j) + \frac{1}{2}(lb^x_i+lb^x_j - lb^y_i - lb^y_j)\sigma_{i,j} \\
    &\geq \min\{lb^x_i+lb^x_j,lb^y_i+lb^y_j\},
\end{align*}
which yields exactly the same constraint set that defines the level-2 combinatorial bound.
\endproof

\end{appendix}

\end{document}